\documentclass{article}

\usepackage[T2A]{fontenc}
\usepackage{color}
\usepackage[utf8]{inputenc} 
\usepackage[english]{babel}
\usepackage{amsmath,amssymb,amsthm}
\usepackage[left=3cm,right=3cm,
top=3cm,bottom=2cm,bindingoffset=1cm]{geometry}
\usepackage{subcaption}
\usepackage{pb-diagram}
\usepackage{graphicx}

\usepackage[utf8]{inputenc} 
\usepackage{amssymb}

\newtheorem{prop}{Proposition}[section]
\newtheorem{lemm}[prop]{Lemma}
\newtheorem{thrm}[prop]{Theorem}
\newtheorem{rmrk}[prop]{Remark}
\newtheorem{defn}[prop]{Definition}
\newtheorem{corl}[prop]{Corollary}
\newcommand{\divg}{{\mathrm{div}\medspace}}
\newcommand{\dom}{{\mathrm{dom}\medspace}}
\newcommand{\curl}{{\mathrm{curl}\medspace}}

\title{Hodge-Laplacian Eigenvalues on Surfaces with Boundary}
\author{Muravyev Mikhail}
\date{}
\begin{document}
	\maketitle
	
	\begin{abstract}
	Recently Rohleder proposed a new variational approach to an inequality between the Neumann and Dirichlet eigenvalues in the simply connected planar case using the language of classical vector analysis. Writing his approach in terms of differential forms permits to generalize these results to a much broader context. The spectrum of the absolute boundary problem for the Hodge-Laplacian on a Riemannian manifold with boundary is presented as a union of the spectra of the absolute boundary problem on the spaces of closed and co-exact forms. An inequality for the eigenvalues of the absolute boundary problem for the Hodge-Laplacian and the Dirichlet boundary problem for the Laplace-Beltrami operator in the Euclidean case is obtained using this presentation. The Rohleder's results are obtained as corollaries of a more general theorem.
	\end{abstract}
	
	\section{Introduction}
	
	\ \ \ \ \ \ Let us consider the following classical question in Spectral Theory in Riemannian Geometry. Given a compact Riemannian manifold with a nonempty boundary, how many Neumann eigenvalues are less then the $k$-th Dirichlet eigenvalue? This question is related, for example, to the the investigation of maxima and minima of eigenfunctions and the hot spots conjecture \cite{Mariano}.
	
	Consider a compact oriented Riemannian manifold $M$ with a boundary $\partial M.$ In this case, the spectra of the Dirichlet and Neumann boundary problems for the Laplace-Beltrami operator are discrete sets without limit points other than $+\infty$. In particular, any Dirichlet or Neumann eigenvalue has finite multiplicity, see \cite[Theorem 2.1.39]{topics}. Let us denote the ordered eigenvalues of  the Dirichlet boundary problem (Dirichlet eigenvalues) by
	$$0 <\lambda_1 < \lambda_2 \leq \lambda_3 \leq \dots,$$
	and the ordered eigenvalues of the Neumann boundary problem (Neumann eigenvalues) by
	$$0 =\mu_1 < \mu_2 \leq \mu_3 \leq \dots.$$
	The eigenvalues are counted accordingly to their multiplicities for both boundary problems, i.e. each eigenvalue is written exactly as many times as the dimension of the corresponding eigenspace. Classical variational descriptions of eigenvalues imply immediately that
	$$\forall k \in \mathbb{N} \quad \mu_k < \lambda_k.$$
	However, a stronger statement, 
	$$\forall k \in \mathbb{N} \quad \mu_{k+1} \leq \lambda_k,$$
	is true for $M \subset \mathbb{R}^n$. It was proven by Friedlander in paper \cite{Friedlander}. Later, Filonov proved in paper \cite{Filonov} that under the same condition there is a strict inequality,
	$$\forall k \in \mathbb{N} \quad \mu_{k+1} < \lambda_k.$$
	For convex Euclidean domains even stronger inequalities are known. If $M \subset \mathbb{R}^2$ is a convex domain with  $C^2$-smooth boundary, then the inequality
	\begin{equation}
		\label{Payne}
		\forall  k \in \mathbb{N} \quad \mu_{k+2} < \lambda_k
	\end{equation}
	holds, it was proven by Payne in \cite{Payne}. Levine and Weinberger showed in \cite{LW} that for a convex $M \subset \mathbb{R}^n$ whose boundary is $C^2$-smooth with H{\"o}lder continuous second derivatives the following estimate holds,
	$$\forall  k \in \mathbb{N} \quad \mu_{k+n} < \lambda_k.$$
	As remarked in \cite{LW}, it can be concluded by approximation argument that the inequality
	$$\forall  k \in \mathbb{N} \quad \mu_{k+n} \leq \lambda_k,$$
	holds for any convex bounded domain $M \subset \mathbb{R}^n.$
	
	The index of the Neumann eigenvalue in the inequality above cannot be greater than $k+n$ without reducing the generality of the result. For instance, for a flat disk one has $\mu_4 > \lambda_1.$ However, conditions on $M$ can still be weakened. For example, Rohleder in his recent work \cite{Rohleder} proved that the convexity condition can be weakened to simply connectedness in case $\dim{M} = 2$. Rohleder's work is interesting not only for its progress in the problem described above, but also because his proof is based on a joint variational description for the Neumann and Dirichlet eigenvalues. 
	
	\begin{defn}
		Consider two nondecreasing sequences $A = \{a_i\}_{i=1}^{+\infty}$ and $B = \{b_i\}_{i=1}^{+\infty}$ such that neither $A$ nor $B$ have limit points but ${+\infty}$. We define \textit{an ordered disjoint union of} $A$ and $B$ as a nondecreasing sequence $C = \{c_i\}_{i=1}^{+\infty}$ such that for any $r \in \mathbb{R}$ we have 
		$$\#\{c_i = r| i\in \mathbb{N}\} = \#\{a_i = r| i\in \mathbb{N}\} + \#\{b_i = r| i\in \mathbb{N}\}.$$
		Here $\# X$ is the cardinality of the set $X.$ We denote such a union by the square cup,
		$$C:= A \sqcup B.$$
	\end{defn}
	Let us denote the ordered disjoint union of the Dirichlet spectrum and the nonzero part of the Neumann spectrum as follows,
	$$\{\eta_k\}_{k=1}^{+\infty} = \{\lambda_k\}_{k=1}^{+\infty} \sqcup \{\mu_k\}_{k=2}^{+\infty}.$$
	Then the following variational description holds.
	
	\begin{thrm}[Rohleder, \protect{\cite[Theorem 4.1]{Rohleder}}]
		\label{rohl_1} 
		Let $M\subset\mathbb{R}^2$ be a bounded simply connected planar domain. Then
		
		$$\eta_k = \min\limits_{\substack{W_k\subset\mathcal{H}_a,\\ \dim(W_k) = k}}\max\limits_{v\in W_k}\frac{\int\limits_{M}\left(\divg^2  (v)+ |\omega (v)|^2\right)}{\int\limits_{M}|v|^2},$$
		where $\omega(v) =\partial_1 v_2 - \partial_2 v_1$  is the vorticity of a vector field $v$ and $\mathcal{H}_{a}$ is the space of all vector fields $v\in L^2(M)^2$ 
		such that $\divg\! v,$ $\omega(v) \in L^2(M)$ and $v|_{\partial M}$  is tangent to the boundary.
	\end{thrm}
	
	Another recent Rohleder's result concerns the problem of finding the eigenvalues of the $\curl \curl$ operator. Let $M \subset \mathbb{R}^3$ be a compact Euclidean domain with Lipschitz boundary $\partial M,$ and consider a vector field $u \in \Gamma(TM)$. The problem is given by the following system,
	\begin{equation}
		\begin{cases}
			\label{curl_problem}
			\curl\curl u = \theta u \ \ \ \ \text{in} \  M,\\
			\divg u = 0 \ \ \ \ \ \ \ \ \ \ \ \ \text{in} \ M,\\
			u \times \nu = 0 \ \ \ \ \ \ \ \ \ \ \ \text{on}\  \partial M.
		\end{cases}
	\end{equation} 
	Here $\nu$ is an exterior normal vector.\\
	
	Let $\{\theta_k\}_{k=1}^{+\infty}$ be the spectrum of problem (\ref{curl_problem}). Then we have the following theorem.
	
	\begin{thrm}[Rohleder, \protect{\cite[Theorem 3.1]{Rohleder_2}}] 
		\label{rohl_2}
		Let $M\subset\mathbb{R}^3$ be a bounded connected Euclidean domain. Then for any $k \in \mathbb{N}$ it is true that
		$$\theta_{2k+1} \leqslant \lambda_k.$$		
	\end{thrm}
	
	Both of the Rohleder`s approaches use the similar technique, but none of them is universal, and each time Rohleder constructs an auxiliary operator anew. 
	
	In this paper we move from vector fields used by Rohleder to differential forms and obtain Theorems \ref{main_result} and \ref{comparisson}. In Section 4 we show that Theorems \ref{rohl_1} and \ref{rohl_2} are special cases of Theorem 4.1. Now let us introduce Theorems \ref{main_result} and \ref{comparisson}.
	
	Consider the absolute boundary condition on the space of all differential p-forms. This problem is given by the following system   
	$$
	\begin{cases}
		\Delta_p\omega = \alpha\omega,\\
		\textbf{n}\omega = 0,\\
		\textbf{n}d\omega = 0,
	\end{cases}
	$$ 
	where $\Delta_p$ is the Hodge-Laplace operator and $\textbf{n}\omega = 0$ means that form $\omega$ is ``tangential''  to the boundary $\partial M$. For a precise definition of $\textbf{n}\omega$ see Section 2.

	Let $D$ be some set of differential $p$-forms. We denote the spectrum of the absolute boundary problem restricted to $L^2D$ by $spec_A(D).$

	\begin{thrm}[\protect{Theorem \ref{main_result}}]
		Let $M$ be a compact oriented Riemannian manifold with a boundary $\partial M$. Then $$spec_A(\Omega_N^p(M)) = spec_A(cE_N^p(M)) \sqcup spec_A(C_N^p(M)).$$ Moreover, if $M$ is simply connected, then $$spec_A(cE_N^p(M)) = spec_A(C_N^{p+1}(M)).$$
		Here $\Omega^p_N(M)$ is a space of all $p$-forms $\omega$ such that $\textbf{n}\omega = 0$, $C^p_N(M)$ is a space of all closed $p$-forms $\omega$ such that $\textbf{n}\omega = 0$ and $cE^p_N(M)$ is a space of all co-exact $p$-forms $\omega$ such that $\textbf{n}\omega = 0.$
	\end{thrm}

	\begin{thrm}[\protect{Theorem \ref{comparisson}}]
		Let $M\subset\mathbb{R}^n$ be a compact Euclidean domain with a smooth boundary $\partial M.$
		Let $\{\theta_k\}^{+\infty}_{k = 1} := spec_A(C^{n-1}_N(M))$. Let $\{\lambda_k\}^{+\infty}_{k = 1}$ be the Dirichlet spectrum of $M$. Then for any $m \in \mathbb{N}$ the following inequality holds,
		$$
		\theta_{(n-1)m + 1} \leqslant \lambda_m.
		$$
	\end{thrm}
		
	\begin{rmrk}
		This paper is an extended English version of the author's M.Sc. thesis \cite{Muravyev} published on-line in Russian in the summer 2024.
	\end{rmrk}

	\begin{rmrk}
		There is a recent preprint \cite{sweden} by M. Fries, M. Goffeng and G. Miranda. The results of \cite{sweden} highly intersects with ones of the present paper. However the results of the present paper are obtained independently and by using another technique.
	\end{rmrk}
	
	\begin{rmrk}
		There is another recent preprint \cite{third_place} by B. Hua, F. Munch, and H. Zhang appeared on the arxiv the same day as the first version of the present paper. The Rohleder's result is generalized in \cite{third_place} from Euclidean domains to hyperbolic surfaces. The proposed generalization of an index shift in (\ref{Payne}) uses the first Betti number. Unlike \cite{third_place}, the present paper is devoted to generalizing to the case of an arbitrary dimension of $M$.
	\end{rmrk}
	
	\section{Basic facts and notations}
	
	\ \ \ \ \ \ Instead of a Euclidean domain let us consider an $n$-dimensional compact oriented Riemannian manifold $(M, g)$ with a boundary $\partial M$. Let $\Omega^p(M)$ be the space of all smooth differential $p$-forms on $M$. We denote the differential, the co-differential and Hodge star operators by $d, \delta$ and $\ast$ respectively. It is worth to recall, that
	$$\delta = (-1)^{n(p+1) + 1}\ast d \ast,$$
	$$\ast \ast = (-1)^{p(p-n)}.$$
	
	Also, $\Omega^p(M)$ can be equipped with the $L^2$ inner product
	$$\langle\omega, \eta\rangle = \int\limits_M \omega\wedge\ast\eta.$$
	In the case, when $\partial M = \varnothing$ operators $d$ and $\delta$ are adjoint, but in general we have the following variation of Stokes` formula. 
	\begin{equation}
		\label{stokes}
		\langle d\omega, \eta\rangle = \int\limits_{M} d\omega \wedge \ast \eta = \int\limits_{M} \omega \wedge \ast \delta\eta + \int\limits_{\partial M} \omega \wedge \ast \eta = \langle\omega, \delta\eta\rangle + \int\limits_{\partial M} \omega \wedge \ast \eta.
	\end{equation}
	
	We write $C^p(M)$ and $cC^p(M)$ for the subspaces of $\Omega^p(M)$ consisting of closed and co-closed $p$-forms, i.e.
	$$C^p(M) := \{\eta \in \Omega^p| d\eta = 0\},$$
	$$cC^p(M) := \{\eta \in \Omega^p| \delta\eta = 0\},$$
	and we write $E^p(M)$ and $cE^p(M)$ for the subspaces of exact and co-exact $p$-forms, i.e.
	$$E^p(M) := \{\eta \in \Omega^p| \eta = d\omega, \omega \in \Omega^{p-1}(M)\},$$
	$$cE^p(M) := \{\eta \in \Omega^p| \eta = \delta\omega, \omega \in \Omega^{p-1}(M)\}.$$
	
	Intersections of these spaces are denoted by a juxtaposition of symbols. For example, $CcC^p(M) = C^p(M) \cap cC^p(M).$ 
	
	The main operator we use is the \textit{Hodge-Laplace} operator,
	$$\Delta_p :\Omega^p(M) \to \Omega^p(M),$$
	$$\Delta_p := \delta d + d\delta.$$
	
	It is well-known that the Hodge-Laplace operator, acting on 0-forms, coincides with the Laplace-Beltrami operator, i.e.
	$$\Delta_0 = \Delta.$$
	
	In the case when $M$ has a nonempty boundary, the Hodge-Laplace operator acting on whole $\Omega^p(M)$ is not self-adjoint, so we need some boundary condition for $p$-forms. We define the map 
	$$\textbf{t}: \Gamma\left(\Lambda^p\left(T^\ast M\right)|_{\partial M}\right)\longrightarrow\Gamma\left(\Lambda^p\left(T^\ast M\right)|_{\partial M}\right)$$
	by a formula 
	$$\textbf{t}\eta(X_1,...,X_p) = \eta(X_1^{||},...,X_p^{||}),\quad \forall X_1, ... X_p \in \Gamma\left(TM|_{\partial M}\right),$$
	where $X = X^{||} + X^{\perp}$ is the decomposition of the vector field $X$ along $\partial M$. Then we define the map $\textbf{n}$ as
	
	$$\textbf{n}: \Gamma\left(\Lambda^p\left(T^\ast M\right)|_{\partial M}\right)\longrightarrow\Gamma\left(\Lambda^p\left(T^\ast M\right)|_{\partial M}\right),$$
	$$\textbf{n}\eta = \eta|_{\partial M} - \textbf{t}\eta.$$

	\begin{rmrk}
		Let $i: \partial M \rightarrow M$ be the inclusion map of the boundary. Abusing notation, $\textbf{t}\eta$ is sometimes identified with the pullback $i^{\ast}$ of the form $\eta$. 
		Note that $\textbf{t}\eta$ belongs to the space $\Gamma\left(\Lambda^p\left(T^\ast M\right)|_{\partial M}\right)$ but $i^{\ast}\eta$ belongs to the space $\Omega^p(\partial M)$, and since they are two different spaces we cannot have $\textbf{t}\eta$ equal to $i^{\ast}\eta$. The case of $\textbf{n}$ and $\iota_\nu$ is similar. The 
		following two propositions clarify why these identifications still make sense.
	\end{rmrk} 
	\begin{prop}[\protect{\cite[Proposition 5.1]{Olle}}] Let $(M, g)$ be a Riemannian manifold with boundary, let $i: \partial M \rightarrow M$ be the inclusion map of the boundary, and let $\eta \in \Omega^p(M).$ Then
		$$\textbf{t}\eta = 0 \Leftrightarrow i^\ast\eta = 0.$$
	\end{prop}
	
	\begin{prop}
		Let $(M, g)$ be a Riemannian manifold with boundary, let $\nu$ be the outward normal vector to the boundary, and let $\eta \in \Omega^p(M).$ Then
		$$\textbf{n}\eta = 0 \Leftrightarrow \iota_{\nu}\eta = 0.$$
	\end{prop}
	$\triangleleft$ 
	If $\textbf{n}\eta = 0$ then $\eta|_{\partial M} = \textbf{t}\eta$, so
	$$\iota_\nu\eta( X_1, X_2, \dots) = \eta(\nu, X_1, X_2, \dots)|_{\partial M} = \textbf{t}\eta(\nu, X_1, X_2, \dots) =$$
	$$=\eta(0, X_1^{||}, X_2^{||}, \dots) = 0.$$
	Note that $X^\perp = f\cdot \nu$, where $f$ is a function. So, due to linearity,  $$\textbf{n}\eta(X_1^{||} + X_1^{\perp}, X_2^{||} + X_2^{\perp}, \dots) = (\eta(X_1^{||} + X_1^{\perp}, X_2^{||} + X_2^{\perp}, \dots) - \eta(X_1^{||}, X_2^{||}, \dots))|_{\partial M}$$
	is a sum of terms of form $f_j \cdot \iota_\nu\eta(\nu, \nu, ..., X^{||}_{i_1}, X^{||}_{i_2}, ...)$. Therefore, if $\iota_{\nu}\eta = 0$ we have $\textbf{n}\eta = 0$. $\triangleright$
	
	We say that a $p$-form $\eta$ satisfies the $\textit{Dirichlet}$ boundary condition if $\textbf{t}\eta = 0.$ If $\textbf{n}\eta = 0,$ we say that $\eta$ satisfies the $\textit{Neumann}$ boundary condition. The corresponding subspaces of $\Omega^p(M)$ are denoted by $\Omega^p_D(M)$ and $\Omega^p_N(M),$ i.e.
	
	$$\Omega_D^p\left(M\right) = \left\{\eta \in \Omega^p\left(M\right)|\textbf{t}\eta = 0 \right\},$$
	$$\Omega_N^p\left(M\right) = \left\{\eta \in \Omega^p\left(M\right)|\textbf{n}\eta = 0 \right\}.$$
	
	Note that the Dirichlet and Neumann boundary conditions for $p$-forms do not generalize the Dirichlet and Neumann boundary conditions for functions. Indeed, $\Omega_N^0 = C^\infty(M)$ is just the space of all $C^\infty$-smooth functions. Moreover, restriction from $\Omega^p$ to $\Omega^p_N$ or $\Omega^p_D$ is still not enough to obtain a self-adjoint operator. However, these conditions still have some good properties.
	
	\begin{prop}[\protect{\cite[Proposition 5.2]{Olle}}] Let $(M, g)$ be a Riemannian manifold with boundary and $\eta \in \Omega^p(M).$ Then
		$$\ast(\textbf{n}\eta) = \textbf{t}(\ast\eta),$$
		$$\ast(\textbf{t}\eta) = \textbf{n}(\ast\eta).$$
	\end{prop}

	\begin{lemm}
		\label{adj}
		Consider $\omega \in \Omega_N^p(M)$ and $\eta\in\Omega_N^{p+1}(M)$ then 
		$$\langle d\omega, \eta\rangle = \langle \omega, \delta\eta\rangle.$$ 
	\end{lemm}
	$\triangleleft$
	Indeed, we have $\int\limits_{\partial M} \omega \wedge \ast \eta = 0$ due to $\textbf{t}\ast\eta = 0.$ Applying (\ref{stokes}) we get
	$$\langle d\omega, \eta\rangle = \langle \omega, \delta\eta\rangle + \int\limits_{\partial M} \omega \wedge \ast \eta = \langle \omega, \delta\eta\rangle.$$ 
	$\triangleright$
	
	\begin{prop}[\protect{\cite[Proposition 5.4]{Olle}}] Let $(M, g)$ be a Riemannian manifold with boundary. Then \\
		\begin{itemize}
			\item the differential $d: \Omega^p(M) \rightarrow \Omega^{p+1}(M)$ preserves the Dirichlet boundary condition,
			\item the co-differential $\delta: \Omega^p(M) \rightarrow \Omega^{p-1}(M)$ preserves the Neumann boundary condition. 
		\end{itemize}
	\end{prop}
	
	Some more notation for different subspaces of $\Omega^p(M),$
	
	$$E^p_D = d\left(\Omega_D^{p-1}\left(M\right)\right) = \{\eta \in \Omega^p\left(M\right)| \eta = d\zeta, \zeta \in \Omega_D^{p-1}\left(M\right)\};$$
	$$cE^p_N = \delta\left(\Omega_N^{p+1}\left(M\right)\right) = \{\eta \in \Omega^p\left(M\right)| \eta = \delta\zeta, \zeta \in \Omega_N^{p+1}\left(M\right)\};$$
	$$C^p_D\left(M\right) = C^p\left(M\right) \cap \Omega^p_D\left(M\right);$$
	$$cC^p_D\left(M\right) = cC^p\left(M\right) \cap \Omega^p_D\left(M\right);$$
	$$C^p_N\left(M\right) = C^p\left(M\right) \cap \Omega^p_N\left(M\right);$$
	$$cC^p_N\left(M\right) = cC^p\left(M\right) \cap \Omega^p_N\left(M\right).$$

	Unfortunately, $\Omega^p(M)$ is not complete with respect to this inner product.
	\begin{defn}
		Consider a subset $A \subseteq \Omega^p(M)$. We denote $L^2$-completion of $A$ by $L^2A.$ 
	\end{defn}
	
	Now we can present some decomposition results.
	
	\begin{thrm}[Hodge-Morrey decomposition, \protect{\cite[Theorem 2.4.2]{Schwarz}}] The space $L^2\Omega^p\left(M\right)$ decomposes into the $L^2$-orthogonal direct sum
		$$L^2\Omega^p\left(M\right) = L^2E_D^p\left(M\right) \oplus L^2CcC^p\left(M\right) \oplus L^2cE_N^p\left(M\right).$$
	\end{thrm}

	\begin{lemm}
		\label{lemma_C_perp_cE}
		The spaces $L^2C^p_N(M)$ and $L^2cE^p_N(M)$ are orthogonal to each other.
	\end{lemm}
	$\triangleleft$
	Due to the Cauchy-Bunyakovsky inequality, it suffices for us to check that $C^p_N(M) \perp cE_N^p(M)$. Let $\omega \in C^p_N(M),$ and $\eta \in cE_N^p(M).$ Then $d\omega = 0$ and there is such a $\xi \in \Omega_N^{p+1}$ that $\eta = \delta\xi.$ Then
	$$
		\langle\omega, \eta\rangle = \langle\omega, \delta\xi\rangle = \int\limits_{M} \omega \wedge\ast\delta\xi = \int\limits_{M} d\omega\wedge\ast\xi - \int\limits_{\partial M} \omega\wedge\ast\xi = 0 + 0 = 0.
	$$
	So, $C^p_N(M)$ is orthogonal to $cE_N^p(M)$, and ,therefore, $L^2C^p_N(M)$ is orthogonal to $L^2cE^p_N(M).$ $\triangleright$
	
	\begin{lemm}
		\label{decompose_1}
		The space $L^2\Omega^p_N(M)$ can be decomposed into the following orthogonal direct sum,
		$$
			L^2\Omega^p_N(M) = L^2C^p_N(M) \oplus L^2cE^p_N(M).
		$$
	\end{lemm}
	$\triangleleft$
	Let $\pi$ be the projection map
	$$
		\pi: L^2\Omega^p(M) \to L^2\Omega^p_N(M).
	$$
	Recall the Hodge-Morrey decomposition for  $L^2\Omega^p(M)$,
	$$
		L^2\Omega^p(M) = L^2E_D^p(M) \oplus L^2CcC^p(M) \oplus L^2cE_N^p(M).
	$$
	We have
	$$
		\left(L^2E_D^p(M) \oplus L^2CcC^p(M)\right) \subset L^2C^p(M),
	$$
	so
	$$
		\pi\left(L^2E_D^p(M) \oplus L^2CcC^p(M)\right) \subset L^2C_N^p(M).
	$$
	Due to Lemma \ref{lemma_C_perp_cE}, we have $L^2C_N^p(M) \perp L^2cE_N^p(M)$, so 
	$$
		L^2C_N^p(M) \subset \left(L^2E_D^p(M) \oplus L^2CcC^p(M)\right)
	$$
	and $\pi(L^2C_N^p(M)) = L^2C_N^p(M)$, which means that 
	$$
		\pi\left(L^2E_D^p(M) \oplus L^2CcC^p(M)\right) = L^2C_N^p(M)
	$$
	and 
	$$
		\Omega^p_N(M) = \pi\left(L^2E_D^p(M) \oplus L^2CcC^p(M) \oplus L^2cE_N^p(M)\right) = L^2C_N^p(M) \oplus L^2cE_N^p(M).
	$$
	$\triangleright$
	
	\section{Boundary problems}
	\ \ \ \ \ \ We make a great use of the following two boundary condition for the Hodge Laplacian acting on $L^2\Omega^p_N(M)$ and  $L^2\Omega^p_0(M) := L^2(\Omega_N^p(M) \cap \Omega_D^p(M))$.  
	\subsection{Absolute boundary condition}
	
	\ \ \ \ \ \ The boundary problem
	$$
	\begin{cases}
		\Delta_p\omega = \alpha\omega,\\
		\textbf{n}\omega = 0,\\
		\textbf{n}d\omega = 0,
	\end{cases}
	$$ 
	for the Hodge-Laplace operator and the corresponding boundary condition are called \textit{absolute}. For reference see \cite{Savo}, \cite[Sections 5.8 - 5.9]{Taylor}.

	We will call eigenvalues and eigenforms of the system above the \textit{absolute} eigenvalues and the \textit{absolute} eigenforms. Also, we will denote the space of all smooth $p$-forms with absolute boundary condition as follows,
	$$\Omega_A^p(M) := \{\omega | \omega \in \Omega^p(M), \textbf{n}\omega = 0, \textbf{n}d\omega = 0\}.$$
	
	\begin{prop}[\protect{\cite[Proposition 5.9.7]{Taylor}}] \label{abs} Eigenvalues of the absolute boundary problem form a discrete set and have no accumulation points but ${+\infty}.$ Moreover, the corresponding eigenforms are analytic and can be chosen to form an orthonormal basis of $L^2\Omega^p_N(M)$. 
	\end{prop}
	
	We will denote the ordered eigenvalues of the absolute boundary problem counted according to their multiplicities by
	
	$$\alpha_1 \leq \alpha_2 \leq \dots \leq \alpha_k \leq \dots$$
	
	\begin{prop}[\protect{\cite[Section 5.9]{Taylor}}]
		The Hodge Laplacian with the absolute boundary have the following quadratic form, 
		$$q[v,u] := \langle \Delta_p u, v\rangle_{L^2} = \langle d u, d v\rangle_{L^2} + \langle \delta u, \delta v\rangle_{L^2}.$$
	\end{prop}
	
	\begin{corl} 
		\label{var_1}
		There is the following variational description for $\{\alpha_k\}_{k=1}^{+\infty}$,
		
		\begin{equation*}
			\alpha_k = \min\limits_{\substack{W_k\subset L^2\Omega^p_N(M),\\ \dim{W_k} = k}}\max\limits_{v\in W_k}\frac{\int\limits_{M}\left(dv\wedge \ast dv + \delta v\wedge \ast \delta v\right)}{\int\limits_{M}v\wedge\ast v}.
		\end{equation*}
	\end{corl}	
	
	\begin{lemm}
		\label{temp_1}
		The differential operator $d$ maps  $\Omega^p_A(M)$ to $\Omega^{p+1}_A(M).$
	\end{lemm}
	$\triangleleft$ Let $\omega$ be an absolute $p$-form, then $\textbf{n}\omega = \textbf{n}d\omega = 0.$ Moreover, $d d \omega = 0$. Therefore,
	\begin{equation*}
		\begin{cases}
			\textbf{n}d\omega = 0,\\
			\textbf{n}dd\omega = 0.
		\end{cases}
	\end{equation*}
	In the other words, $d\omega$ belongs to $\Omega_A^{p+1}(M)$.
	$\triangleright$
	
	\begin{defn}
		We define the spaces of closed and exact absolute as follows,
		$$C_A^p(M) = \{\omega| \omega\in\Omega_A^p(M), d\omega = 0\},$$
		$$E_A^p(M) = Im_d(\Omega_A^{p-1}).$$
	\end{defn} 
	
	\begin{defn}
		We define the absolute cohomology as follows,
		$$H^p_A(M) := C_A^p(M)/E_A^p(M).$$
	\end{defn}
	The quotion space in the definition above make sense due to Lemma \ref{temp_1} .
	\begin{defn}
		We define the space of absolute harmonic $p$-forms as follows:
		$$\mathcal{H}^p_A(M) := \{\omega \in \Omega_A^p(M)|d \omega = \delta\omega = 0\}.$$
	\end{defn}

	\begin{prop}[\protect{\cite[Proposition 5.9.10]{Taylor}}] \label{noharm} Consider the space of absolute harmonic $p$-forms $\mathcal{H}^p_A(M)$ and the absolute cohomology $H^p_A(M)$, then 
		$$H^p_A(M) \cong \mathcal{H}^p_A(M).$$
		Hence, in the simply connected case
		$$\mathcal{H}^p_A(M) = \{0\}.$$
	\end{prop}
	
	Now let us move to the main result of this paper.
	
	\begin{defn}
		Let $D$ be a subset of $\Omega_N^k.$ We denote the spectrum of the absolute boundary problem restricted to $L^2D$ by $spec_A(D).$
	\end{defn} 
	
	\begin{thrm}
		\label{main_result}
		Let $M$ be a compact oriented Riemannian manifold with a boundary $\partial M$. Then $$spec_A(\Omega_N^p(M)) = spec_A(cE_N^p(M)) \sqcup spec_A(C_N^p(M)).$$ Moreover, if $M$ is simply connected, then $$spec_A(cE_N^p(M)) = spec_A(C_N^{p+1}(M)).$$
	\end{thrm}
	$\triangleleft$ 
	Due to Proposition \ref{abs} there is an orthonormal basis $A^p$ of the space $L^2\Omega_N^p(M),$ where $A^p$ consists only of the absolute eigenforms. Let us show that $A^p$ can be chosen in a form $A^p = A^p_{cE} \cup A^p_C,$ where $A^p_{cE}$ is an orthonormal basis of the space $L^2cE_N^p(M)$ and $A^p_{C}$ is an orthonormal basis of the space $L^2C_N^p(M).$ Here the union of sequences means that we are considering a new sequence, where the elements of the first sequence are in even places, and the elements of the second are in odd places 
	
	Let $\omega$ be an absolute $p$-eigenform with an eigenvalue $\alpha$. Due to Lemma \ref{decompose_1} there are unique $\eta \in L^2C^p_N(M)$ and $\xi \in L^2cE_N^p(M)$ such that $\omega = \eta + \xi.$ Note that $\delta \xi = d \eta = 0$, so there are two following equalities,
	$$
		\Delta_p(\xi) = \delta d \xi,\quad
		\Delta_p(\eta) = d\delta \eta.
	$$
	Let us compute $\Delta_p$ on $\omega,$
	$$
		\alpha(\eta + \xi) = \alpha \omega = \Delta_p(\omega) = \Delta_p(\eta + \xi) = d\delta\eta + \delta d\xi. 
	$$
	Due to the uniqueness of $\xi$ and $\eta$ we have the following equalities,
	$$
		\Delta_p(\xi) = \delta d \xi = \alpha \xi,\quad
		\Delta_p(\eta) = d\delta \eta = \alpha \eta.
	$$
	Therefore, $\eta$ and $\xi$ are absolute eigenforms. The author leaves the process of further reconstruction of the basis $A^p$ to the reader.
	
	As long as we can chose $A^p$ in the form $A^p = A^p_C \cup A^p_{cE}$ the following decomposition holds, $$spec_A(\Omega_N^p(M)) = spec_A(cE_N^p(M)) \sqcup spec_A(C_N^p(M)).$$
	
	Now let us prove, that if $M$ is simply connected, then
	$$
		Spec_A(C^{p+1}_N(M)) = Spec_A(cE^p_N(M)).
	$$
	First of all, let us show that $Spec_A(C^{p+1}_N(M)) \subset Spec_A(cE^p_N(M))$. Here $\subset$ means that if the sequence $Spec_A(C^{p+1}_N(M))$ contains $m$ copies of an element $\alpha$ then $Spec_A(cE^p_N(M))$ contains at least $m$ copies of an element $\alpha$.\\
	Let $\omega_1\perp\omega_2$ be two orthogonal absolute $(p+1)$-eigenforms. Let $\Delta_{p+1}\omega_1 = \alpha_1\omega_1$ and $\Delta_{p+1}\omega_2 = \alpha_2\omega_2$. Let us show that $\delta\omega_1$ and $\delta\omega_2$ are nonzero absolute $p$-eigenforms and $\delta\omega_1\perp\delta\omega_2$, $\Delta_p\delta\omega_1 = \alpha_1\delta\omega_1$ and $\Delta_p\delta\omega_2 = \alpha_2\delta\omega_2$.
	
	Forms $\delta\omega_1$ and $\delta\omega_2$ are nonzero, because elsewhere either $\omega_1$ or $\omega_2$ is harmonic, and there is a contradiction with Lemma \ref{noharm}.
	
	Now we compute $\Delta_p\delta\omega_1,$
	\begin{align*}
		\Delta_p\delta\omega_1 = (\delta d + d\delta)\delta\omega_1 = \delta d \delta\omega_1 = (\delta\delta d + \delta d \delta)\omega_1 = \delta \Delta_{p+1}\omega_1 = \alpha_1 \delta \omega_1.
	\end{align*}
	
	Let us show that $\delta\omega_1$ satisfy the absolute boundary condition. The condition $\textbf{n}\delta\omega_1 = 0$ holds because $\textbf{n}\omega_1 = 0$ and the operator $\delta$ preserves the Neumann boundary condition for forms. Moreover, $\textbf{n}d\delta\omega_1 = \textbf{n}\alpha_1\omega_1 = 0.$
	
	So we can conclude that $\delta\omega_1$ and $\delta\omega_2$ are the absolute $p$-eigenforms. Now let us show that $\delta\omega_1 \perp \delta\omega_2.$ Due to Lemma \ref{adj} we have
	
	\begin{align*}
		\langle\delta\omega_1, \delta\omega_2\rangle = \langle d\delta\omega_1, \omega_2\rangle = \alpha_1 \langle\omega_1, \omega_2\rangle = 0.
	\end{align*}   
	
	Now let $C_\alpha^{p+1}$ be an eigenspace of the operator $\Delta_{p+1}$ with a domain $C_N^{p+1}(M)$ and $cE_\alpha^{p}$ be an eigenspace of the operator $\Delta_{p}$ with a domain $cE_N^{p}(M).$ Both of the eigenspaces correspond to the eigenvalue $\alpha$. Note, that in such notations the orthonormal basis of $C_\alpha^{p+1}$ which consists of the absolute $(p+1)$-eigenforms generates an orthonormal basis of $cE_\alpha^p$, which consists of the absolute $p$-eigenforms. So we have
	\begin{align*}
		spec_A(C_\alpha^{p+1})\subset spec_A(cE_\alpha^p).
	\end{align*}   
	Hence,
	\begin{align*}
		spec_A(C_N^{p+1}(M))\subset spec_A(cE_N^p(M)).
	\end{align*} 
	
	Now we prove that $spec_A(C_N^{p+1}(M)) \supset spec_A(cE_N^p(M)).$
	
	Let us show that for any co-exact absolute $p$-eigenform $\xi\in cE^p_N(M)$ with an eigenvalue $\alpha$ there is a such closed absolute $(p+1)$-eigenform $\omega\in C^{p+1}_N(M)$ with the eigenvalue $\alpha$ that $\delta \omega = \xi.$
	
	The form $\xi$ is co-exact, so there is $p+1$-form  $\eta \in \Omega_N^{p+1}(M)$ such that $\delta\eta = \xi.$ Due to Lemma \ref{decompose_1} there are such $p+1$-forms $\omega_1 \in L^2C_N^{p+1}(M)$ and $\omega_2 \in L^2cE_N^{p+1}(M)$ that $\eta = \omega_1 + \omega_2.$ Since $\delta\omega_2 = 0,$ we know that $\xi = \delta\omega_1.$
	
	Let us show that $\omega_1$ is an absolute $(p+1)$-eigenform with the eigenvalue $\alpha$,
	
	$$
		\alpha\delta\omega_1 = \Delta_p \delta \omega_1 = (d \delta \delta  + \delta d \delta) \omega_1 = (\delta \delta d + \delta d \delta) \omega_1 = \delta \Delta_{p+1}\omega_1.
	$$
	Therefore, $\Delta_{p+1}\omega_1 = \alpha\omega_1 + \chi$, where $\chi \in L^2cE_N^{p+1}(M)$.
	
	Note that $d\omega_1 = 0$, hence, $\Delta_{p+1}\omega_1 = d\delta\omega_1 \in L^2C_N^{p+1}(M)$, and $\chi = 0.$ Therefore, $\omega_1$ is an absolute $p$-eigenform with the eigenvalue $\alpha$ and $\delta\omega_1 = \xi.$ So, in a way similar to the proof that  $spec_A(C_N^{p+1}(M))\subset spec_A(cE_N^p(M))$, we can conclude that $spec_A(C_N^{p+1}(M))\supset spec_A(cE_N^p(M))$. Hence,
	$$
		spec_A(C_N^{p+1}(M)) = spec_A(cE_N^p(M)).
	$$
	$\triangleright$
	
	In the end of this section we propose the following lemma, which will be useful in obtaining the Rohleder`s results. 
	
	\begin{lemm}
		\label{aux_lemm_1}
		Let $M$ be a compact oriented Riemannian manifold of dimension $n$ with a boundary $\partial M.$ Then $spec_A(\Omega_N^n(M)) = spec_A(C^n_N(M)) = \{\lambda_k\}_{k=1}^{+\infty}.$ Here $\{\lambda_k\}_{k=1}^{+\infty}$ is the Dirichlet spectrum on $M.$
	\end{lemm}
	$\triangleleft$
	Note that any $n$-form on $M$ is closed, so $C_N^n(M) = \Omega_N^n(M).$\\
	
	Any $n$-form in $\Omega_N^n(M)$ can be writen as $fdV\!\!ol$. For such a forms the absolute boundary condition reduces to $\textbf{n}(fdV\!\!ol) = 0$, or, in the other terms, $\textbf{t}f = 0$. Condition $\textbf{t}f = 0$ on a functions is equal to $f|_{\partial M} = 0$, so 
	$$
		spec_A(\Omega_N^n(M)) = \{\lambda_k\}_{k=1}^{+\infty}.
	$$
	$\triangleright$
	
	\subsection{Second boundary condition}	
	\ \ \ \ \ \ Consider the space $\Omega^p_0(M) := \Omega_N^p(M) \cap \Omega_D^p(M).$ One can consider the following boundary condition and the corresponding boundary problem on $\Omega^p_0(M)\subset\Omega^p_N(M),$ 	
	
	\begin{align}
		\label{true}
		\begin{cases}
			\Delta_p\omega = \alpha\omega,\\
			\textbf{n}\omega = 0,\\
			\textbf{t}\omega = 0.
		\end{cases}
	\end{align}
	
	Due to the reasons described below, throughout this paper we will call (\ref{true}) the \textit{true Dirichlet} condition  (in order to separate it from the usual Dirichlet boundary condition.) The corresponding eigenvalues and eigenforms will be called  \textit{true Dirichlet} too.
	
	\begin{prop}
		\label{sym}
		Operator $\Delta_p$ with domain $\Omega^p_0(M)$ is symmetric and positive.
	\end{prop}
	$\triangleleft$ For any $\omega, \eta \in \Omega^p_0(M)$
	we have 
	$$\langle\Delta_p \omega, \eta\rangle = \int\limits_M \Delta\omega \wedge \ast \eta = \int\limits_M d \delta\omega \wedge \ast \eta + \int\limits_M \delta d \omega \wedge \ast \eta = \int\limits_M d \delta\omega \wedge \ast \eta + \int\limits_M \eta \wedge \ast \delta d \omega = $$
	
	$$= \int\limits_M \delta\omega \wedge \ast \delta\eta + \int\limits_M d\eta \wedge \ast d \omega + \int\limits_{\partial M}\delta\omega \wedge \ast \eta + \int\limits_{\partial M} \eta \wedge \ast d\omega = $$
	
	$$=\langle\delta\omega, \delta\eta\rangle + \langle d\eta, d \omega\rangle + \int\limits_{\partial M}\delta\omega \wedge \ast \eta - \int\limits_{\partial M} \eta \wedge \ast d\omega.$$
	Note that $\eta|_{\partial M} = \ast\eta|_{\partial M} = 0.$ Therefore,
	
	$$\langle\Delta_p \omega, \eta\rangle = \langle\delta\omega, \delta\eta\rangle + \langle d\eta, d \omega\rangle.$$
	
	Thus, we have shown that $q[\omega,\eta] := \langle\Delta_p \omega, \eta\rangle$ is a quadratic form, so, for any $\omega, \eta \in \Omega_0^p$ we have $\langle\Delta_p \omega, \eta\rangle = \langle \omega,\Delta_p \eta\rangle$ and $\langle\Delta_p \omega, \omega\rangle \geq 0.$
	$\triangleright$
	
	\begin{corl}[Friedrichs extension \protect{\cite[Theorem X.23]{Reed}}] In the notation of Proposition \ref{sym}, $q$ is a closable quadratic form and its closure $\bar{q}$ is a quadratic form of a unique self-adjoint operator $\bar{\Delta}_p$. Operator $\bar{\Delta}_p$ is a positive extension of $\Delta_p$, and the lower bound of its spectrum is the lower bound of $q$.    
	\end{corl}
	
	Since $\bar{\Delta}_p$ is self-adjoint, the eigenvalues of the true Dirichlet boundary problem form a discrete set and have no accumulation points but ${+\infty}.$ We will denote the ordered eigenvalues of the true Dirichlet boundary problem counted according to their multiplicities by
	
	$$\rho_1 \leq \rho_2 \leq \dots \leq \rho_k \leq \dots$$
	
	Also, we have the following variational description for $\{\rho_k\}_{k=1}^{+\infty}$,
	
	\begin{align}
		\label{var_pho}
		\rho_k = \min\limits_{\substack{W_k\subset \dom(\bar{\Delta}_p),\\ \dim{W_k} = k}}\max\limits_{v\in W_k}\frac{\int\limits_{M}\left(dv\wedge \ast dv + \delta v\wedge \ast \delta v\right)}{\int\limits_{M}v\wedge\ast v}.
	\end{align}
	
	In the Euclidean case it is quite easy to build an orthonormal basis of $\Omega_0^p(M)$ consisting of true Dirichlet eigenforms. 
	
	\begin{thrm}
		Consider a bounded Euclidean domain $M\subset \mathbb{R}^n$ and $p$-forms $\zeta^i_k := u_k dx^{i_1}\wedge\dots\wedge dx^{i_p}$, where  $i_1 < i_2 < \dots < i_p$ and ${u_k}_{k=1}^{+\infty}$ is the orthonormal basis of $C^{\infty}_0$ consisting of the Dirichlet eigenfunctions. Then the following two statements hold.
		\begin{itemize}
			\item Any $p$-form $\zeta^i_k$ is the a Dirichlet eigenform with an eigenvalue $\lambda_k.$
			\item Set of all $p$-forms $\zeta^i_k$ form an orthonormal basis in the space $\Omega^p_0$.
		\end{itemize} 
	\end{thrm}
	
	$\triangleleft$
	Again, the completeness follows from the completeness of $\{u_k\}_{k=1}^{+\infty}$ in $\Omega_N^0(M).$ Moreover, $\zeta^i_k$ satisfies the true Dirichlet boundary condition, since $u_k|_{\partial M} = 0.$
	
	Let us show the orthogonality. If for $\zeta^i_k$ and $\zeta^j_m$ one has $i\neq j$, then $$dx^{i_1}\wedge\dots\wedge dx^{i_p} \wedge \ast dx^{j_1}\wedge\dots\wedge dx^{j_p} = 0.$$
	Hence,
	$$\langle\zeta^i_k, \zeta^j_m \rangle = \int\limits_M 0 = 0.$$
	For $\zeta^i_k$ and $\zeta^i_m$ with same $i$, we have
	$$\langle\zeta^i_k, \zeta^i_m \rangle = \int\limits_M u_k u_m dV\!\!ol = \delta_m^k.$$
	
	Now, let us show, that $\zeta^i_k$ is an eigenform. Let
	
	$$\ast dx^{i_1}\wedge\dots\wedge dx^{i_p} := dx^{j_1}\wedge\dots\wedge dx^{j_{n-p}}.$$
	Then,
	$$\ast dx^{j_1}\wedge\dots\wedge dx^{j_{n-p}} = (-1)^{p(n-p)} dx^{i_1}\wedge\dots\wedge dx^{i_p};$$
	
	$$\ast dx^{j_b}\wedge dx^{i_1}\wedge\dots\wedge dx^{i_p} = (-1)^{p+b-1} dx^{j_1}\wedge\dots\wedge \hat{dx^{j_b}}\wedge\dots\wedge dx^{j_{n-p}};$$
	
	$$\ast dx^{i_a}\wedge dx^{j_1}\wedge\dots\wedge dx^{j_{n-p}} = (-1)^{(p+1)(n-p) + a -1} dx^{i_1}\wedge\dots\wedge \hat{dx^{i_a}}\wedge\dots\wedge dx^{i_p}.$$
	Therefore, 
	
	$$\delta d (\zeta^i_k) = (-1)^{n(p+2)+1} \ast d \ast\left(\sum\limits_{b=1}^{n-p}\frac{\partial u_k}{\partial x^{j_b}}dx^{j_b}\wedge dx^{i_1}\wedge\dots\wedge dx^{i_p}\right) = $$
	$$= (-1)^{np + p + b} \ast d \left(\sum\limits_{b=1}^{n-p}\frac{\partial u_k}{\partial x^{j_b}}dx^{j_1}\wedge\dots\wedge \hat{dx^{j_b}}\wedge\dots\wedge dx^{j_{n-p}}\right)=$$
	$$= (-1)^{np + p + b} \ast \left(\sum\limits_{a=1}^{p}\sum\limits_{b=1}^{n-p}\frac{\partial^2 u_k}{\partial x^{j_b}\partial x^{i_a}}dx^{i_a}\wedge dx^{j_1}\wedge\dots\wedge \hat{dx^{j_b}}\wedge\dots\wedge dx^{j_{n-p}}\right) + $$
	$$+ (-1)^{np + p + 1} \ast \left(\sum\limits_{b=1}^{n-p}\frac{\partial^2 u_k}{\partial x^{j_b}\partial x^{j_b}} dx^{j_1}\wedge\dots\wedge dx^{j_{n-p}}\right) = $$
	$$= (-1)^{a + 1}  \left(\sum\limits_{a=1}^{p}\sum\limits_{b=1}^{n-p}\frac{\partial^2 u_k}{\partial x^{j_b}\partial x^{i_a}}dx^{j_b}\wedge dx^{i_1}\wedge\dots\wedge \hat{dx^{i_a}}\wedge\dots\wedge dx^{i_p}\right) - $$
	$$ - \left(\sum\limits_{b=1}^{n-p}\frac{\partial^2 u_k}{\partial x^{j_b}\partial x^{j_b}}dx^{i_1}\wedge\dots\wedge  dx^{i_p}\right),$$
	and
	
	$$d\delta(\zeta^i_k) = (-1)^{n(p+1) +1} d \ast d \left(u_k dx^{j_1}\wedge\dots\wedge dx^{j_{n-p}}\right) = $$
	$$ = (-1)^{n(p+1) +1} d \ast \left(\sum\limits_{a = 1}^{p}\frac{\partial u_k}{\partial x^{i_a}} dx^{i_a}\wedge dx^{j_1}\wedge\dots\wedge dx^{j_{n-p}}\right)=$$
	$$ = (-1)^{a} d \left(\sum\limits_{a = 1}^{p}\frac{\partial u_k}{\partial x^{i_a}} dx^{i_1}\wedge\dots\wedge \hat{dx^{i_a}}\wedge\dots\wedge dx^{i_p}\right)=$$
	$$ = (-1)^{a}\left(\sum\limits_{b = 1}^{n - p}\sum\limits_{a = 1}^{p}\frac{\partial^2 u_k}{\partial x^{i_a}\partial x^{j_b}} dx^{j_b}\wedge dx^{i_1}\wedge\dots\wedge \hat{dx^{i_a}}\wedge\dots\wedge dx^{i_p}\right) -$$
	$$ - \left(\sum\limits_{a = 1}^{p}\frac{\partial^2 u_k}{\partial x^{i_a}\partial x^{i_a}} dx^{i_1}\wedge\dots\wedge  dx^{i_p}\right).$$
	Finally, 
	$$\Delta_p(\zeta^i_k) = (\delta d + d\delta)\zeta^i_k =  -\sum\limits_{a = 1}^{n}\frac{\partial^2 u_k}{\partial x^{a}\partial x^{a}}dx^{i_1}\wedge\dots\wedge  dx^{i_p} = \Delta_0(u_k)dx^{i_1}\wedge\dots\wedge  dx^{i_p} = \lambda_k \zeta^i_k.$$
	$\triangleright$
	
	Since $\{\zeta_k^i\} \subset \Omega_0^p(M)$ is also a basis in  $\Omega_0^p(M),$ it is a basis of $\dom\bar{\Delta}_p$. Therefore, in the variational principle (\ref{var_pho}) instead of $\dom\bar{\Delta}_p$ we can write $\Omega_0^p(M),$ 
	
	$$\rho_k = \min\limits_{\substack{W_k\subset\Omega^p_0(M),\\ \dim{W_k} = k}}\max\limits_{v\in W_k}\frac{\int\limits_{M}\left(dv\wedge \ast dv + \delta v\wedge \ast \delta v\right)}{\int\limits_{M}v\wedge\ast v}.$$
	
	Also, from Theorem 3.2, we know, that $\{\rho_k\}_{k=1}^{+\infty}$ is $\left(^p_n\right)$ copies of $\{\lambda_k\}_{k=1}^{+\infty}$, so we have the following statement.
	
	\begin{prop}
		\label{var_2}
		Consider a bounded Euclidean domain $M\subset\mathbb{R}^n$ and let $\{\rho_k\}_{k=1}^{+\infty}$ be equal to the ordered disjoint union $\bigsqcup\limits_{i = 1}^{\left(^p_n\right)}\{\lambda_k\}_{k=1}^{+\infty},$ then 
		
		$$\rho_k = \min\limits_{\substack{W_k\subset\Omega^p_0(M),\\ \dim{W_k} = k}}\max\limits_{v\in W_k}\frac{\int\limits_{M}\left(dv\wedge \ast dv + \delta v\wedge \ast \delta v\right)}{\int\limits_{M}v\wedge\ast v}.$$
		
	\end{prop}
	
	\section{Variational principles comparison}
	We begin our comparison of the two boundary problems with a rather general theorem. 
	\begin{thrm}
		\label{comparisson}
		Let $M\subset\mathbb{R}^n$ be a compact Euclidean domain with a smooth boundary $\partial M.$
		Let $\{\theta_k\}^{+\infty}_{k = 1} := spec_A(C^{n-1}_N(M))$. Let $\{\lambda_k\}^{+\infty}_{k = 1}$ be the Dirichlet spectrum of $M$. Then for any $m \in \mathbb{N}$ the following inequality holds,
		$$
			\theta_{(n-1)m + 1} \leqslant \lambda_m.
		$$
	\end{thrm}
	$\triangleleft$ 
	At first, let $m\in \mathbb{N}$ be such an integer that $\lambda_m < \lambda_{m+1}$.
	
	Let denote the whole absolute spectrum on $\Omega^n_N(M)$ as 
	$$\{\alpha_k\}^{+\infty}_{k = 1} := spec_A(\Omega_N^p(M)).$$
	Theorem \ref{main_result} and Lemma \ref{aux_lemm_1} imply
	$$
		\{\theta_k\}^{+\infty}_{k = 1} \sqcup \{\lambda_k\}^{+\infty}_{k = 1} = \{\alpha_k\}^{+\infty}_{k = 1}.
	$$
	Corollary \ref{var_1} provides the following variational description,
	\begin{align}
		\label{var_3}
		\alpha_k = \min\limits_{\substack{W_k\subset \Omega^n_N(M),\\ \dim{W_k} = k}}\max\limits_{\substack{v\in W_k,\\v \neq 0 }}\frac{\int\limits_{M}\left(dv\wedge \ast dv + \delta v\wedge \ast \delta v\right)}{\int\limits_{M}v\wedge\ast v}.
	\end{align}	
	On the other hand, due to Proposition \ref{var_2} there is a variational description for $n$ copies of the Dirichlet spectrum. Let $\{\rho_k\}^{+\infty}_{k = 1}:= \bigsqcup\limits_{i = 1}^{n}\{\lambda_k\}_{k=1}^{+\infty},$ then 	
	\begin{align}
		\label{var_4}
		\rho_k = \min\limits_{\substack{W_k\subset \Omega^p_0(M),\\ \dim{W_k} = k}}\max\limits_{\substack{v\in W_k,\\v \neq 0 }}\frac{\int\limits_{M}\left(dv\wedge \ast dv + \delta v\wedge \ast \delta v\right)}{\int\limits_{M}v\wedge\ast v}.
	\end{align}	
	Let $\{u_k\}_{k=1}^{mn}$ be first $mn$ eigenforms of degree $n-1$, corresponding to the true Dirichlet boundary condition and to the first $mn$ elements of $\{\rho_k\}^{+\infty}_{k = 1}$. Then each $u_k$ belongs to $\Omega_0^{n-1}(M)\subset\Omega_N^{n-1}(M)$, so we can choose $span(\{u_k\}_{k=1}^{mn})$ as a subspace $W_{nm}$ in  (\ref{var_3}) and obtain
	$$
		\alpha_{mn} \leq \lambda_m,
	$$
	and, therefore,
	\begin{align}
		\label{ineq_1}
		\theta_{m(n-1)} \leq \lambda_m,
	\end{align}
	Now let us show, that lower index of the left hand side of (\ref{ineq_1}) can be increased by 1. Let $v$ be an absolute $(n-1)$-eigenform with the eigenvalue $\lambda_m.$ First of all, note that $v$ is linearly independent from $span(\{u_k\}_{k=1}^{mn})$ because otherwise $v$ would satisfy the following three boundary conditions,
	$$
		\begin{cases}
			\textbf{n}v = 0,\\
			\textbf{t}v = 0,\\
			\textbf{n}dv = 0,
		\end{cases}
	$$
	and since operator $d$ preserves the Dirichlet boundary condition on forms, the form $v$ would satisfy the fourth boundary condition, 
	$$
		\textbf{t}dv = 0.
	$$
	In the other words, $v$ would satisfy 
	$$
		\begin{cases}
			v|_{\partial M} = 0,\\
			dv|_{\partial M} = 0,
		\end{cases}
	$$	
	and then, due to the fact that $v$ is analytic, we get $v = 0.$
	Secondly, due to the fact that $v$ is the absolute eigenform with the eigenvalue $\lambda_m$, we know that for any $(n-1)$-form $\omega \in \Omega_N^{n-1}(M)$ it is true that
	
	\begin{equation*}
		\langle d v, d\omega \rangle + \langle \delta v, \delta\omega \rangle = q[v, \omega] = \langle \Delta_{n-1} v, \omega \rangle = \lambda_m \langle v, \omega \rangle 
	\end{equation*}
	and hence $q[\omega, \eta] \leqslant \lambda_m \langle \omega, \eta \rangle $ for any $\omega, \eta \in span(\{u_k\}_{k=1}^{mn}\cup\{v\}).$
	
	Therefore, we can chose $W_{mn+1}$ for (\ref{var_3}) as $span(\{u_k\}_{k=1}^{mn}\cup\{v\})$ and obtain  $$\alpha_{mn+1} \leqslant \lambda_m.$$ 
	As long as $\lambda_m < \lambda_{m+1}$, we get  $$\theta_{m(n-1) + 1} \leq \lambda_m.$$
	
	Now let $\lambda_m = \lambda_{m+1}.$ Let $l$ be the first natural number greater than $m$ such that $\lambda_{l} < \lambda_{l+1},$ then
	$$
		\theta_{m(n-1) + 1} \leqslant \theta_{l(n-1) + 1} \leqslant \lambda_{l} = \lambda_m, 
	$$
	hence, 
	$$
		\theta_{m(n-1) + 1} \leqslant\lambda_m. 
	$$
	$\triangleright$
	
	Now let us show, how to obtain the Rohleder's results as particular cases of Theorem \ref{comparisson} in dimensions 2 and 3.
	\begin{corl}[\protect{\cite[Theorem 4.1]{Rohleder}}]
		Let $M \subset \mathbb{R}^2$ be a compact simply connected planar domain of dimension 2 with a smooth boundary $\partial M$. Let $\{\lambda_k\}_{k = 1}^{+\infty}$ be the Dirichlet spectrum on $M$ and $\{\mu_k\}_{k = 1}^{+\infty}$ be the Neumann spectrum on $M$. Then for all $k \in \mathbb{N}$ the following inequality holds,
		$$
		\mu_{k+2} \leqslant \lambda_k.
		$$
	\end{corl}
	$\triangleleft$ In the case, when $dim M = 2$, we have the following equality of spaces,
	$$\Omega_N^{n-2}(M) = \Omega_N^0(M) = C^{\infty}(M).$$
	Here $C^{\infty}(M)$ corresponds to the space of all $C^\infty$ smooth functions on $M$.	The equality holds since for all $f \in C^{\infty}(M)$ it is true that $\textbf{n}f = 0.$\\
	We know that for all $f \in C^{\infty}(M)$
	$$df = 0 \Leftrightarrow f = const,$$
	hence, the space of closed $0$-forms consists only of the constant functions, 
	$$ C^0_N(M) = C^0(M) = \{c | c\in \mathbb{R}\}.$$
	Therefore,
	$$cE^0_N(M) = \{c | c\in \mathbb{R}\}^{\perp}.$$
	
	On functions the absolute boundary condition reduces to the Neumann boundary condition, because $\textbf{n}f = 0$ holds for all $f \in C^{\infty}(M)$ and 
	$$\textbf{n}df = 0 \Leftrightarrow \frac{\partial f}{\partial \nu}|_{\partial M} = 0.$$
	Hence, $spec_A(\Omega_N^0(M)) = \{\mu_k\}^{+\infty}_{k=1}$ and $$spec_A(cE^0_N(M)) = \{\mu_k\}^{+\infty}_{k=1} \setminus spec_A(\{c|c\in\mathbb{R}\}) = \{\mu_k\}^{+\infty}_{k=1} \setminus \{0\} = \{\mu_k\}^{+\infty}_{k=2}.$$
	Theorem \ref{comparisson} and Theorem \ref{main_result} imply
	$$\mu_{k+2} \leqslant \lambda_k.$$
	$\triangleright$
	
	\begin{corl}[\protect{\cite[Theorem 3.1]{Rohleder_2}}]
		Let $M \subset \mathbb{R}^3$ be a compact Euclidean domain of dimension 3 with a smooth boundary $\partial M$. Let $\{\lambda_k\}_{k = 1}^{+\infty}$ be the Dirichlet spectrum on $M$ and $\{\theta_k\}_{k = 1}^{+\infty}$ be the spectrum of the boundary problem given by the system
		\begin{equation}
			\label{curlcurl}
			\begin{cases}
				v \in \Gamma(TM),\\
				\curl \curl v = \theta v,\\
				\divg v = 0,\\
				v \times \nu = 0.
			\end{cases}
		\end{equation}
		Then for all $k \in \mathbb{N}$ the following inequality holds,
		$$
		\theta_{2k+1} \leqslant \lambda_k.
		$$
	\end{corl} 
	$\triangleleft$
	Consider the canonic isomorphism $\varphi$ between the spaces $\Gamma(TM)$ and $\Omega^1(M)$ induced by the Euclidean metric. For $\omega\in \Omega^1(M)$ and $v \in \Gamma(TM)$ the equality $\varphi(v) = \omega$ holds if and only if for all $u \in \Gamma(TM)$ we have $\omega(u) = \langle u, v\rangle$. In Cartesian coordinates, the isomorphism has the form
	$$\varphi\left(f_x\frac{\partial}{\partial x} + f_y\frac{\partial}{\partial y} + f_z\frac{\partial}{\partial z}\right) = f_xdx + f_ydy + f_zdz.$$
	
	Let $$\omega := f_xdx + f_ydy + f_zdz$$
	and $$v := f_x\frac{\partial}{\partial x} + f_y\frac{\partial}{\partial y} + f_z\frac{\partial}{\partial z},$$
	then
	\begin{equation*}
		\divg v =\frac{\partial f_x}{\partial x} + \frac{\partial f_y}{\partial y} + \frac{\partial f_z}{\partial z};
	\end{equation*}

	\begin{equation*}
		\delta \omega = (-1)^{3\times2+1}\ast d \ast \omega = 
		- (\ast d(f_x dy\wedge dz + f_y dz\wedge dx + f_z dx\wedge dy)) = 
	\end{equation*}
	\begin{equation*}
		= -\left(\ast\left(\frac{\partial f_x}{\partial x}dx\wedge dy \wedge dz + \frac{\partial f_y}{\partial y}dy\wedge dz \wedge dx + \frac{\partial f_z}{\partial z}dz\wedge dx \wedge dy\right)\right) = -\left(\frac{\partial f_x}{\partial x} + \frac{\partial f_y}{\partial y} + \frac{\partial f_z}{\partial z}\right);
	\end{equation*}

	\begin{equation*}
		\curl v = \left(\frac{\partial f_z}{\partial y} - \frac{\partial f_y}{\partial z}\right)\frac{\partial}{\partial x} + \left(\frac{\partial f_x}{\partial z} - \frac{\partial f_z}{\partial x}\right)\frac{\partial}{\partial y} + \left(\frac{\partial f_y}{\partial x} - \frac{\partial f_x}{\partial y}\right)\frac{\partial}{\partial z};
	\end{equation*}

	\begin{equation*}
		\ast dv = \ast\left(
		\frac{\partial f_x}{\partial y}dy \wedge dx + \frac{\partial f_x}{\partial z}dz \wedge dx + \frac{\partial f_y}{\partial x}dx \wedge dy + \frac{\partial f_y}{\partial z}dz \wedge dy + \frac{\partial f_z}{\partial x}dx \wedge dz + \frac{\partial f_z}{\partial y}dy \wedge dz\right) = 
	\end{equation*}
	\begin{equation*}
		= \left(\frac{\partial f_z}{\partial y} - \frac{\partial f_y}{\partial z}\right) dx + \left(\frac{\partial f_x}{\partial z} - \frac{\partial f_z}{\partial y}\right) dy + \left(\frac{\partial f_y}{\partial x} - \frac{\partial f_x}{\partial y}\right) dz.
	\end{equation*}
	Therefore, we have
	\begin{enumerate}
		\item $\divg v = 0 \Leftrightarrow \delta \omega = 0$;
		\item $v \times \nu = 0 \Leftrightarrow \textbf{t} \omega = 0$;
		\item $\varphi(\curl v) = \ast d\omega.$
	\end{enumerate}
	We obtain from 1. and 3. that
	$$\varphi(\curl \curl v) = \delta d\omega = \Delta_1\omega.$$ 
	Therefore, problem (\ref{curlcurl}) on $v$ is equal to the following problem on $\omega$. 
	\begin{equation*}
		\begin{cases}
			\omega \in \Omega^1(M),\\
			\Delta_1\omega = \theta \omega,\\
			\delta \omega = 0,\\
			\textbf{t}\omega = 0.
		\end{cases}
	\end{equation*}
	The problem above is equal to its dual problem
	\begin{equation*}
		\begin{cases}
			\omega \in \Omega^1(M),\\
			\Delta_2\ast\omega = \theta \ast\omega,\\
			d\ast\omega = 0,\\
			\textbf{n}\ast\omega = 0.
		\end{cases}
	\end{equation*}
	Obviously, we have 
	$$d\ast\omega = 0 \Rightarrow \textbf{n}d\ast\omega = 0,$$
	Hence, the system (\ref{curlcurl}) is equal to the system
	\begin{equation*}
		\begin{cases}
			\omega \in \Omega^1(M),\\
			\Delta_2\ast\omega = \theta \ast\omega,\\
			d\ast\omega = 0,\\
			\textbf{n}d\ast\omega = 0,\\
			\textbf{n}\ast\omega = 0,
		\end{cases}
	\end{equation*}
	And it is the absolute boundary condition, restricted on $C^2_N(M),$ therefore, $spec_A(C_N^2(M)) = \{\theta_k\}_{k = 1}^{+\infty}$ and Theorem \ref{comparisson} implies
	$$\theta_{2k+1} \leqslant \lambda_k.$$
	$\triangleright$\\

	\textit{Acknowledgments.} The author would like to express many thanks to  Mikhail Karpukhin, Alexei Penskoi and Iosif Polterovich for the great support they provided during the writing of this paper.
	

\end{document}